\documentclass[reqno,11pt]{amsart}

\usepackage{amsmath}
\usepackage{amssymb,esint}
\usepackage{graphicx}

\newtheorem{theorem}{Theorem}[section]
\newtheorem{corollary}[theorem]{Corollary}
\newtheorem{lemma}[theorem]{Lemma}
\newtheorem{proposition}[theorem]{Proposition}
\theoremstyle{definition}
\newtheorem{definition}[theorem]{Definition}

\newtheorem{remark}[theorem]{Remark}

\numberwithin{equation}{section}

\usepackage{enumitem}
\setenumerate{label={\rm (\alph{*})}}

\usepackage{hyperref,enumitem,mathrsfs}
\usepackage[rose,novargreek]{paper_diening}
\usepackage{macros}

\DeclareMathAlphabet\mathbfcal{OMS}{cmsy}{b}{n}
\DeclareMathAlphabet\mathbfscr{OMS}{mdugm}{b}{n}

\providecommand{\curl}{\mathrm{curl}\,}
\providecommand{\divergence}{\mathrm{div}}
\providecommand{\Div}{\divergence\,}
\newcommand{\N}{\mathbf{N}}
\newcommand{\R}{\mathbb{R}}

\title[Weak solutions for the Rotational Smagorinsky model]{On the existence of
  weak solutions for a family of unsteady  rotational Smagorinsky models}

\author[L.C. Berselli]{Luigi C. Berselli}

\address[L.C.~Berselli]{Dipartimento di Matematica, Universit\`a di
  Pisa, Via F. Buonarroti 1/c, I-56127 Pisa,  Italy}
\email{{\tt luigi.carlo.berselli@unipi.it}}

\author[A.~Kaltenbach]{Alex Kaltenbach}
\address[A.~Kaltenbach and M.~R\accent23u\v{z}i\v{c}ka]
{Institute of Applied Mathematics,
  Albert-Ludwigs-University Freiburg,
  Ernst-Zermelo-Stra{\ss} 1, 79104 Freiburg,}
\email{\tt alex.kaltenbach@mathematik.uni-freiburg.de}
\email{\tt rose@mathematik.uni-freiburg.de}
\author[R. Lewandowski]{Roger Lewandowski}
\address[R. Lewandowski]{IRMAR, UMR 6625, Universit\'e Rennes 1,
   Campus Beaulieu, 35042 Rennes cedex FRANCE}
\email{\tt   Roger.Lewandowski@univ-rennes1.fr}
\author[M. R\accent23u\v{z}i\v{c}ka]{Michael R\accent23u\v{z}i\v{c}ka}



\begin{document}

\begin{abstract}
  In this paper we show that the rotational Smagorinsky model for turbulent flows, can be
  put, for a wide range of parameters in the setting of Bochner pseudo-monotone evolution
  equations. This allows to prove existence of weak solutions a) identifying a proper
  functional setting in weighted spaces and b) checking some easily verifiable
  assumptions, at fixed time. We also will discuss the critical role of the exponents
  present in the model (power of the distance function and power of the curl) for what
  concerns the application of the theory of pseudo-monotone
  operators. 
\end{abstract}

\maketitle


\section{Introduction}
In this paper we introduce the unsteady general
rotational Smagorinsky model for incompressible turbulence
\begin{equation}
  \label{eq:Baldwin--Lomax-alpha}
\begin{aligned} 
      \partial_{t}\vm+\vortm\times\vm +\curl\big(C_\alpha \ell^\alpha|\vortm|\vortm \big)
      +\nabla \overline{q}&=\bff&& \quad\mbox{in
        $(0,T)\times\Omega$,}
      \\
      \vortm&=\curl\vm&&\quad \mbox{in        $(0,T)\times\Omega$,}
      \\
      \Div \vm&=0&& \quad\mbox{in $(0,T)\times\Omega$,}
      \\
      \vm&=\mathbf{0}&&\quad \mbox{on $(0,T)\times\partial\Omega$,}
      \\
      \vm(0)&=\overline{ \vv_{0} }&& \quad\mbox{in $\Omega$,}
    \end{aligned}
\end{equation}
where $\Om$ is a smooth bounded domain in $\R^3$, $\ell$ is the
Prandtl mixing length, $\alpha >0$ is a given exponent, $C_\alpha>0$
is a calibration constant, $\vm$ is the mean velocity, $\vortm$ is the
mean vorticity, and $\overline{q}$ is the sum of the Bernoulli
pressure of the fluid and certain potentials such as the turbulent
kinetic energy and others. Here and in the sequel, for each smooth
vector $\bu: \R^{3}\to\R^{3}$ we define as curl the vector
\begin{equation*}
  (\curl\bu)_{i}:=\sum_{j,k=1}^{3}\epsilon_{ijk}\frac{\partial
    u_{j}}{\partial x_{k}},
\end{equation*}
where $\epsilon_{ijk}$ is the Levi--Civita totally anti-symmetric tensor.

Note that in the
equations~\eqref{eq:Baldwin--Lomax-alpha} the linear dissipative term
$-\nu\Delta\vm$ is not present, since we are considering flows at very
high Reynolds number, and then viscous effects are negligible compared
to the Reynolds stresses. Note that the
presence of the linear dissipative term, for any $\nu>0$ will allow
for some simplifications of the
proofs. Nevertheless, to obtain result and estimates independent of
$\nu>0$, a treatment as the one we provide is requested.

According to standard
assumptions (see \eqref{eq:Prandtl_Obukhov} and \eqref{eq:van_driest} in
Section \ref{sec:Rot_Smag_Mod}), we will assume that $\ell$ behaves as
the distance to the boundary. This means that $\ell(\x) \approx d(\x)$ when
$\x \in \Om$ and $\x$ is close to the boundary $\p \Om$ (see \eqref{eq:assumption}
and \eqref{eq:assumption_2} in what follows). As it is common in 
turbulence modeling, we assume that the flow fields are  stochastic
processes, and the bar operator stands for the expectation in the
Reynolds decomposition $\vv = \vm + \vv'$, $\pi = \overline{\pi} +
\pi'$, where $\pi $ denotes the pressure and $\overline {\pi}$ the
mean pressure (see
Section~\ref{sec:reyn_decomp} below, even if other choices are
possible, as for instance denoting by the bar operator the long time-averaging).

The natural value of the parameter $\alpha$ is equal to 2 and this
model is similar to the widely used Smagorinsky model, but the term
$ \div (C_s \ell^2 | D \vm | D \vm)$ is replaced here by
$\curl (C_2 \ell^2 | \vortm | \vortm)$. The equivalence between both
models can be understood for homogeneous isotropic turbulence, by the
equality of the enstrophy $ \overline { | \vorticity |^2 } $ to the
total mean deformation $2 \overline {| D \vv | ^2 }$. The equivalence
can be obtained by a straightforward generalization
of~\cite[Lemma~4.7]{CL2014}. Then, according to this equality,
in~\cite[Section~5.5.1]{CL2014} it is proved that the $-5/3$
Kolmogorov law yields to express the eddy viscosity as
$\nu_{T} = C_2 \ell^2 | \vortm |$. The rotational structure of the
eddy diffusion is a peculiarity of the model which is suitable for
high-speed flows with thin attached boundary-layers. The mathematical
treatment of rotational models is one of the main theoretical
contribution of this paper.

The numerical performance of this model in the steady state case has
been initially tested by Baldwin and Lomax~\cite{BL1978}, so that this
model is also known as the Baldwin--Lomax model. Numerical analysis
foundations also in the statistical non-equilibrium setting can be
found in~\cite{RLZ2019}.

The analytical properties of a steady version of this model have been
recently studied in~\cite{BB2020} in the setting of weighted Sobolev
spaces. Some unsteady versions, with the presence of a dispersive
term --which allows for a more classical treatment-- have been recently
studied in~\cite{Ngu2020,BLN2021}.

The steady version can be treated within the standard theory of
monotone operators, plus a localization argument, while the unsteady
one requires a more delicate argument to deal with the precise choice
of spaces and formulation of the problem. As we will prove, a proper
definition of the functional setting will make
system~\eqref{eq:Baldwin--Lomax-alpha} to fit into the framework of
evolution problems with Bochner pseudo-monotone operators, for which
the theory have been recently developed by two of the authors
in~\cite{KR2021}. The theory developed in~\cite{KR2021} represents an
extension and an adaption to unsteady problems of the classical theory
of pseudo-monotone operators from Br\'ezis~\cite{Bre66}, \cite{Bre1968}, already
described in the classical monograph of Lions~\cite{Lio1969}. Our main
result is the following, which covers all possible positive powers of
the distance function which are strictly smaller than the critical value $\alpha=2$.
\begin{theorem}\label{thm:main-theorem}
  Let us suppose that $\ell(\x)=d(\x,\partial\Omega)$ and let $\alpha\in [0,2)$, $0<T<\infty$,
  $\overline{ \vv_0 } \in L^2_\sigma(\Omega)$, and
  $\bff\in L^{3/2}(0,T;(W^{1,3}_0(\Omega,d^\alpha)^*)$. Then, there
  exists a weak solution to the initial boundary value problem
  \eqref{eq:Baldwin--Lomax-alpha} such that
\begin{equation*}
  \begin{aligned}
    & \vm\in 
    C([0,T];L^2_\sigma(\Omega))\cap
    L^3(0,T;W^{1,3}_{0,\sigma}(\Omega,d^\alpha)),
    \\
    &\text{and for all }t\in[0,T]
    \\
    &\frac{1}{2}\|\vm(t)\|^2+\int_0^t\int_\Omega C_\alpha
    d^\alpha(\bx)|{\vortm }
    (s,\bx)|^3\,d\bx\,ds=\frac{1}{2}\|\overline{ \vv_0}
    \|^2+\int_0^t\langle\bff,\vm\rangle_{W^{1,3}_0(\Omega,d^\alpha)}ds.
      \end{aligned}
\end{equation*}
\end{theorem}
The limitation $\alpha<2$ seems to be intrinsic to the problem due to
the fact that $d^{\alpha}$ is not anymore a Muckenhoupt weight for
$\alpha\geq2$ (cf.~Definition \ref{def:muc}). Hence, for $\alpha\geq2$ most of the analytical
properties may fail, since we cannot ensure that the quantity from the
energy estimates controls the (weighted) full gradient of the
solution. For values of $\alpha$ larger or equal than $2$, even the
weak formulation, the density of smooth functions, and the meaning of
the boundary conditions may fail; the solution of the problem, if
possible, would pass through the introduction of a more general setting,
of very weak solutions.

In the last section we will also consider the existence for a family of problems with
different powers of the vorticity in the turbulent stress tensor, still with the distance
function raised to any exponent smaller than the critical one,
cf. Thm.~\ref{thm:main-theorem-p}.

\bigskip

\textbf{Plan of the paper.} In Section \ref{sec:modeling} we derive
the Rotational Smagorinsky model from a classical turbulence modeling process, in Section \ref{sec:sec3} we define
the notion of Bochner pseudo-monotone operators and we recall the main result
for general evolutionary problems, in Section~\ref{sec:weighted} we
recall the main results on weighted spaces, which will be used to
properly formulate the problem. 
Next in
the final Section~\ref{sec:sec5} we show how the hypotheses apply to problem
\eqref{eq:Baldwin--Lomax-alpha}, for relevant choices of the weight
functions and discuss generalization and critical values of the parameters.

\section{Modeling}\label{sec:modeling}  

\subsection{Reynolds decomposition} \label{sec:reyn_decomp} Let us
consider the Navier--Stokes equations (NSE in the sequel) written with
the convective term 
in the rotational formulation: 
\begin{equation} 
  \label{eq:NS_standard_form3}  
  \begin{aligned}
    \vv_t + \vorticity \times \vv - \nu \Delta \vv + \g \left(\pi + \dfrac{|\vv|^2}{2}
    \right) &= \fv&&\quad \mbox{in
        $(0,T)\times\Omega$,}
    \\
    \vorticity&= \curl\vv&&\quad \mbox{in
        $(0,T)\times\Omega$,}
    \\
    \Div \vv &= 0&&\quad \mbox{in
      $(0,T)\times\Omega$,}
    \\
    \vv 
    &= {\bf 0}&&\quad \mbox{on
      $(0,T)\times\partial\Omega$,}
    \\
    \vv (0) 
    & = \vv_0&&\quad \mbox{in
        $\Omega$,}
  \end{aligned} 
\end{equation}
where $\vv = \vv(t, \x, \omega)$ is the velocity field,
$\pi = \pi(t, \x, \omega)$ the pressure, $\vorticity = \curl \vv$ the
vorticity, $(t, \x) \in \R_+ \times \Om$,
$\omega \in X({\mathcal B}, P)$, where $X({\mathcal B}, P)$ is a given
probability space on the space of initial data.

For instance, if $\bff=\mathbf{0}$ (the argument can be adapted also
to include a smooth
enough external force) it holds
that for each divergence-free element of $\bfv_{0}\in H^{1/2}(\Omega)$
there exists a lower bound $T=T(\| \bfv_0\|_{1/2})>0$  for the life-span of the unique
Fujita--Kato mild solution. Since the life-span can be estimated with the norm of
the initial datum, by fixing $X=\overline{B(0,R)}\subseteq H^{1/2}(\Omega)\cap
\{\nabla\cdot\vv=0\}$ for some $R>0$, then the life-span is bounded
from below by some $T_{X}>0$. This means that for each $\vv_{0}\in X$,
there exists a unique $\vv\in C(0,T_{X};H^{1/2}(\Omega))$ solution of
the NSE.

We then introduce $P$, which is a probability measure on the Borel
sets~of~$X$. More specifically, $P$ can be constructed as limit of
averages of Dirac measures as in \cite{CL2014} or the renormalized
Lebesgue measure constructed from the Borel sets of $X$. The final
result does not depend on the choice of $P$.
Let us denote the expectation with a bar, hence
\begin{equation*}
  \overline {\vv_0} = \int_X \vv_0 \,d P(\vv_0),
\end{equation*}
 and 
 \begin{equation*}
\vm (t, \x) = \int_X \vv(t, \x, \vv_0) \,d P(\vv_0)\qquad \overline{\pi} (t, \x)
= \int_X \pi (t, \x, \vv_0)\, d P(\vv_0).
\end{equation*}
 More generally, for any
field $\Psi = \p_t \vv, \vorticity, \g \vv, \Delta \vv, { | \vv |^2 \over 2}...$,
we can define the statistical mean as 
\begin{equation*}
  \overline \Psi(t,\x) = \int_X \Psi(t, \x, \vv_0) \,d P(\vv_0),
\end{equation*}
and consequently we can perform the usual decomposition of  $\Psi$ as 
\begin{equation*}
\Psi = \overline \Psi +\Psi' ,
\end{equation*}
which is known as the Reynolds decomposition. 
The properties of the statistical averaging process imply (Reynolds rules) that  for all
$\Psi,\Theta\in X$
\begin{equation*}
\p_t \overline \Psi  = \overline {\p_t \Psi}, \qquad \g \overline
\Psi   = \overline {\g \Psi }, \qquad  \overline{\Psi'} = 0,
\qquad\overline{\overline{\Psi}\,\Theta}=\overline{\Psi}\,\overline{\Theta},
\end{equation*}
hence, taking the expectation of the NSE~\eqref{eq:NS_standard_form3} yields
\begin{equation} 
  \label{eq:NS_standard_form3_mean}  
    \begin{aligned}
      \overline{\vv}_t + \overline \vorticity \times \vm + \overline{ \vorticity' \times
        \vv'} - \nu \Delta \vm
      +\g \left (\overline \pi + {|\vm|^2 \over 2} + \overline{ |\vv'
          |^2\over 2} \right ) &= \overline{\fv},
      \\
      \vortm &= \curl\vm, 
      \\
      \Div \overline{\vv} &= 0, 
      \\
      \vm \vert_{\p \Om}
      &= {\bf 0}, 
      \\
      \vm \vert_{t=0}
      & = \overline{\vv_0} .
      \end{aligned}
\end{equation}
The basic closure and modeling problems concern expressing
$\overline{ \vorticity' \times \vv'} $~in~terms of averaged variables.
\subsection{Rotational Reynolds stress} 
When taking the expectation of the NSE with the convective term
written in the usual form, we get the term
$\div (\overline{ \vv' \otimes \vv'} )$. The quantity
$\reyn = \overline{\vv' \otimes \vv'}$ is called the Reynolds stress
and the Boussinesq assumption consists in assuming that 
\begin{equation*}
\reyn = - \nu_{T} D \vm,
\end{equation*}
where $\nu_{T}\geq0$ is an eddy viscosity which remains to be
determined and modeled in terms of $\vm$. If we want to use such a
Boussinesq assumption, we must express the turbulent stress
(which is a vector in the rotational formulation)
$$
\bs:=\overline{ \vorticity' \times \vv'}
$$ 
in terms of derivatives of mean quantities. This is similar to the
approach used when modeling the more standard
Reynolds stress tensor.  We prove in what follows the following theorem
\begin{theorem}
  Assume that $\Om$ is connected and of class $C^1$. Then, there
  exists a vector
  $ {\bf a}^{\text{(R)}}= {\bf a}^{(\text{R})} (t , \x)$
  and a scalar potential  $\Phi = \Phi (t, \x)$ such that
  \begin{equation}
    \label{eq:NS_standard_form4_mean}  
      \begin{aligned}
        \overline{\vv}_t + \vortm \times\vm + \curl {\bf a}^{(\hbox
          {\tiny R})} - \nu \Delta \vm + \g (\overline \pi +  \tfrac 1 2\vert \overline \bfv\vert^2+ k - \Phi)
        &= \overline{\fv},
        \\
        \Div \overline{\vv} &= 0,\\
        \vm \vert_{\p \Om} &= {\bf 0},\\
        \vm \vert_{t=0} & = \overline{\vv_0},
      \end{aligned}
  \end{equation}
  where $k = {1 \over 2} \overline {| \vv '
    |^2}$  is the turbulent kinetic energy. 
\end{theorem}
\begin{proof} Let $ {\bf a}^{(\hbox {\tiny R})}$ and $\Phi$ be given by:
  \begin{equation}
    \begin{aligned}
      \label{eq:gen_Reyn_Stress} {\bf a}^{(\hbox {\tiny R})} (t , \x)
      &= {1 \over 4 \pi} \int_\Om {\curl {\bf s}(t, \x') \over |\x-\x'
        | } d \x' + {1 \over 4 \pi} \int_{\p \Om}{ {\bf s}(t, \x')
        \over | \x-\x'| }\times d\sigma (\x'),
      \\
      \Phi (t, \x) &= {1 \over 4 \pi} \int_\Om {\div \,{\bf s}(t, \x')
        \over |\x-\x' | } d \x' - {1 \over 4 \pi} \int_{\p \Om}{ {\bf
          s}(t, \x') \over | \x-\x'| }\cdot d\sigma (\x').
    \end{aligned}
\end{equation}
Therefore, by the Helmholtz--Hodge theorem, we have the relation 
\begin{equation}
\label{eq:gen_Reyn_Stress_2}  \overline{ \vorticity' \times \vv'} =
\curl   {\bf a}^{(\hbox {\tiny R})}  - \g \Phi.
\end{equation}
 Inserting 
\eqref{eq:gen_Reyn_Stress_2} into \eqref{eq:NS_standard_form3_mean}
gives \eqref{eq:NS_standard_form4_mean}.  
\end{proof} 
The vector ${\bf a}^{(\hbox {\tiny
    R})}$, which is continuously and uniquely determined by formula
\eqref{eq:gen_Reyn_Stress}, is called the rotational Reynolds stress
tensor.
From now on we write $\bar q $ instead of $\overline \pi + \tfrac 1 2\vert \overline \bfv\vert^2+ k - \Phi$.

\subsection{Closure assumption: Rotational Smagorinsky model} \label{sec:Rot_Smag_Mod} 
In order to finish the modeling of turbulent quantities, it remains to
link ${\bf a}^{(\hbox {\tiny R})}$ to the mean vorticity $\overline
\vorticity$. Notice that ${\bf a}^{(\hbox {\tiny
    R})}$ has the dimension of a squared velocity, while $\overline
\vorticity$ those of a frequency. Therefore, adapting the Boussinesq
assumption to this case yields to assume
\begin{equation*}
  {\bf a}^{(\hbox {\tiny R})} = \nu_{T}\, \overline
  \vorticity,
\end{equation*}
in which
$\nu_{T}\geq{0}$ is a quantity with the dimensions of a
viscosity. According to the
$-5/3$ Kolmogorov law and following~\cite[Section~5.5.1]{CL2014}, we
can assume  (for an homogeneous and isotropic flow, in the limit
$\nu\to0$) 
$$
\nu_{T} = \nu_{T} (\ell, | \vortm|),
$$
where $\ell$ is the Prandtl mixing length.
The dimensional analysis of the expression shows that a consistent
expression is  
\begin{equation}
\nu_{T} = C_\omega \ell^2 | \vortm |,
\label{eq:nut}
\end{equation}
with $C_\omega$ a dimensionless constant.
This raises the question of the determination of $\ell$.
In the case of a	 flow over a plate, one finds in Obukhov~\cite{Obu1946} the
following classical  law:
\begin{equation} \label{eq:Prandtl_Obukhov} \ell = \ell (z) = \kappa z,
\end{equation}
where $z\geq0$ is the distance from the plate and
$\kappa$ the von K\'arm\'an constant.  The Van Driest
formula~\cite{Dri1956} defines $\ell$ by:
\begin{equation} 
  \label{eq:van_driest} 
    \ell (z) := \kappa\, z\, (1- \mathrm{e}^{-z/A});
\end{equation} 
here $A$ depends on the oscillations of the plate and on the kinematic viscosity $\nu$,
while $z\geq0$ is again the distance from the plate.

According to these formula, we shall assume throughout the rest of the
paper that the function $\ell:\,\overline
\Omega\to\R^{+}$ is of class
$C^{2}$ and satisfies the two following properties:
    \begin{align}
    \label{eq:assumption} 
      a)\ 
       & \ell (\x) \approx d(\x, \partial\Omega) \qquad \text{for
     } \x  \hbox{ close to } \, \partial\Omega;
      \\
 \label{eq:assumption_2} 
 b)\  &\forall \ K \subset \subset \Om, \,\,  \exists \,  \ell_K>0\quad
      \hbox{s.t.} \quad \ell (\x ) \ge \ell_K >0 \quad \forall \, \x
      \in K,    \end{align} 
where $d(\x,\partial\Omega)$ denotes the distance from the
boundary. In practice, we could have directly assumed
$\ell(\x)=d(\x)$, i.e.,
\begin{equation}
\nu_{T} = C_\omega d^2 | \vortm |.
\label{eq:nut1}
\end{equation}
\subsection{Generalised Rotational Smagorinsky models by
  dimensional analysis}
The analysis of the previous section can be put also in a more general
framework of Large Eddy Simulation (LES) models, looking also at
possible modifications of the parameters present in the expression of
the turbulent (rotational) stress vector. Let
$\ell_0>0$ be a typical length scale of the motion. For instance, in
the case of a flow over a plate, one can take
\begin{equation*}
\ell_0 = {\nu \over v_{*}},
\end{equation*}
where $\nu$ is the kinematic viscosity and $v_{*}$ is the so-called
friction velocity  (cf.~\cite{ABLN2019}).

We consider (modulo introducing an appropriate non-dimensiona\-lization
of the equations) the following operator
\begin{equation}\label{eq:p=3}
  \curl\big(\ell_0^{2-\alpha}\ell ^\alpha|\vortm|\vortm\big)
\end{equation}
with $\alpha \in [0,2]$, which is degenerate at the boundary and for which the natural
treatment is through scales of weighted Banach spaces.

We report some discussion about the relationships between the scaling
of the weight and that of the power of the curl. In the framework of
LES methods we show that even starting with
\begin{equation}\label{eq:p}
  \nu _T= \ell_0^{2-\alpha}\ell ^\alpha|\vortm|^{p-2}
\end{equation}
this determines a link between powers $\alpha$ and $p$.
 Nevertheless, in the
last section we will also point out the limiting behavior of the
exponent $p=3$ present in model~\eqref{eq:Baldwin--Lomax-alpha}, when $p=\alpha+1$.

If one thinks of a flow as composed of eddies of different sizes in
different places, then in a region of large eddies the changes of velocity and
its curl  are both $\mathcal{O}(1)$ of the typical distance. In
a region of smaller eddies the velocity changes over a distance of
$\mathcal{O}$(eddy length scale), so the local deformation is
$\mathcal O$(1/eddy length scale),
cf.~\cite[\S~3.3.2]{BIL2006}. Hence, the rotational Smagorinsky model
introduces a turbulent viscosity $\nu_{T}=(C\delta)^{2}|\vortm |$,
where $\delta$ is the (local) smallest resolved scale, such that
\begin{equation*}
  \nu_{T}=\left\{
    \begin{aligned}
      &\mathcal{O}(\delta^{2})&&\qquad\text{in regions where }|\vortm |=\mathcal{O}(1),
      \\
      &\mathcal{O}(\delta)&&\qquad\text{in the smallest resolved scale where
      }|\vortm |=\mathcal{O}(\delta^{-1}).
    \end{aligned}
\right.
\end{equation*}
By extrapolation, motivated by experiments with central difference
approximations to linear convection diffusion problems, the following
alternate scaling has also been proposed (cf.~\cite{BIL2006} and
Layton~\cite{Lay1996}) $\nu_{T}=(C\delta)^{p-1}|\bfD\vm|^{p-2}$, and
we consider here the rotational counterpart
\begin{equation*}
\nu_{T}=(C\delta)^{p-1}|\vortm|^{p-2},\qquad 1<p<\infty,
\end{equation*}
which resembles general power laws for non-Newtonian fluids. 
The above choice of $\nu_{T}$  satisfies 
\begin{equation*}
  \nu_{T}=\left\{
    \begin{aligned}
      &\mathcal{O}(\delta^{p})&&\qquad\text{in regions where
      }|\vortm|=\mathcal{O}(1), 
      \\
      &\mathcal{O}(\delta)&&\qquad\text{in the smallest resolved scale where
      }|\vortm |=\mathcal{O}(\delta^{-1}),
    \end{aligned}
\right.
\end{equation*}
The justification of the presence of the critical value $p-1$ as power
of the distance function can be done directly by dimensional
arguments~as~in~\cite{BB2020}. In fact, recall that both
$\nabla\overline \bfv$ and $\overline \bfomega$ have dimensions
$T^{-1}$, where $T$~is~a~time, and in \eqref{eq:nut} the turbulent
viscosity $\nu_{T}=d^{2}|\vortm|\sim L^{2}T^{-1}$ (where $L$ is a
length) has the dimensions of a viscosity. This is the only way to
identify (by using just a typical length and the vorticity) a quantity
with the dimensions of a viscosity. Introducing as third parameter as
the friction velocity $v_{*}\sim L T^{-1}$, one can consider more
general combinations.  The outcome is to find a turbulent eddy viscosity of
the following form
\begin{equation*}
  \nu_{T}=v_{*}^{\theta}d^{\alpha}|\vortm|^{p-2},
\end{equation*}
for some constants $\theta,\,\alpha,\,p$. It turns out
(cf.~\cite{BB2020}) that the dimensions of this quantity are
$\nu_{T}\sim L^{\theta+\alpha}T^{2-\theta-p}$, and to respect
dimensions of the viscosity one has to fix
\begin{equation*}
  \theta=3-p\qquad\text{and}\qquad \alpha=p-1.
\end{equation*}
A sound  generalization of the rotational
Smagorinsky  model is then the  one with rotational  stress  
\begin{equation*}
  \mathbf {S}(v_{*},d,\vortm)=C v_{*}^{3-p}d^{p-1}|\vortm |^{p-2} \vortm,
\end{equation*}
and, after re-scaling, one can assume $C v_{*}^{3-p}=1$.  Note that, even for different
values of $p$, the power of the distance is always the critical one (in terms of
analytical properties of the weight functions), since $d^{p-1}\not\in {A}_{p}$,
cf.~Lemma~\ref{lem: Muckenhoupt}.

In the last section we will show that from the point of view of mathematical
properties, the turbulent eddy viscosity 
$$
  \nu _T=d^{p-1}|\vortm|^{p-2},
$$
can be handled in terms of an existence
theory by (pseudo)monotone operators only for $p\geq3$. Hence, the exponent $p=3$ plays
for the weighted rotational operators, the same role that the exponent $p=11/5$ plays for
the usual $p$-NSE with stress tensor $S(\bD\vm)=c|\bD\vm|^{p-2}$.

From now and so far no risk of confusion occurs, we do not write the bar anymore. 

\section{Evolution equations in an abstract setting} 
\label{sec:sec3}
As already claimed in the introduction, a proper setting to the
rotational Smagorinsky model is that of pseudo-monotone evolution problems so
we briefly recall the abstract existence result we will use on the
sequel.

For the convenience of the reader, we recall the following~known~definition.
\begin{definition}\label{2.1}
  Let $X$, \hspace*{-0.15em}$Y$ be Banach spaces. An
  operator~${A\hspace*{-0.1em}:\hspace*{-0.2em}X\hspace*{-0.15em}\rightarrow\hspace*{-0.15em} Y}$~is~called
  \begin{enumerate}
  \item[\textrm{(i)}] \textbf{bounded}, if for all bounded
    $M\hspace*{-0.2em}\subseteq\hspace*{-0.2em} X$, the image $A(M)\hspace*{-0.2em}\subseteq\hspace*{-0.2em} Y$~is~bounded.
  \item[\textrm{(ii)}] \textbf{coercive}, if $Y=X^*$ and
    $\lim\limits_{\|x\|_X \to \infty}\frac{ \langle
      Ax,x\rangle _X}{\|x\|_X}=\infty$.
  \item[\textrm{(iii)}] \textbf{pseudo-monotone}, if $Y\hspace*{-0.2em}=\hspace*{-0.2em}X^*$ and for
    any~sequence~${(x_n)_{n\in\mathbb{N}}\hspace*{-0.2em}\subseteq\hspace*{-0.2em} X}$~from
    \begin{align*}
      x_n\overset{n\rightarrow\infty}{\rightharpoonup}x\text{ in }X
      ,&
      \\
      \limsup_{n\rightarrow\infty}{\langle Ax_n,x_n-x\rangle_X}\leq
      0,&
    \end{align*}
    it follows that
    $\langle Ax,x-y\rangle_X\hspace*{-0.1em}\leq\hspace*{-0.1em} \liminf_{n\rightarrow\infty}{\langle
      Ax_n,x_n-y\rangle_X}$~for~all~${y\hspace*{-0.1em}\in\hspace*{-0.1em} X}$.
  \end{enumerate}
\end{definition}
It is well-known that for each $f\in X^{*}$, the steady problem $Ax=f$~admits~a~solution if $A$ is bounded, coercive and pseudo-monotone,
see~ \cite{Bre66}, \cite{Bre1968}. A typical example of a pseudo-monotone operator is
the sum of a hemi-continuous monotone and a compact operator. Recently, two of the authors in~\cite{KR2021}
developed an abstract framework for evolution problems, by using the
concepts~of~Bochner pseudo-monotone and Bochner coercive operators to
generalize the ideas of \cite[Sec. 2.5]{Lio1969},~\cite{landes-87},~\cite{Hir1},~\cite{Hir2}~and~\cite{Shi}. We want to access this theory
for our concrete example.  Therefore, for the remainder
of this section, we assume that $(V,H,\textrm{id})$ is an
evolution triple, i.e., $V$ is a separable, reflexive Banach space,
$H$ a separable Hilbert space
and $V$ embeds densely into $H$. 
For ${I:=\left(0,T\right)}$, $T\in \left(0,\infty\right)$, and
${p\in \left(1,\infty\right)}$,~we~set
\[
 \mathbfcal{X}:=L^p(I,V)\qquad \text{and}\qquad
 \mathbfcal{Y}:=L^\infty(I,H).
\]
In this framework we have the following notion of a time derivative. 
\begin{definition}\label{2.10}
  A
  function $\mathbf{u}\in\mathbfcal{X}$ has a
  \textbf{generalized time derivative} 
  if there exists a function
  ${\mathbf{w}\in L^{p'}(I,V ^*)}$ such that 
  \begin{align*}
	-\int_I{(\mathbf{u}(s),v)_H\varphi^\prime(s)\,ds}=\int_I{\langle\mathbf{w}(s),v\rangle_{V}\varphi(s)\,ds}
  \end{align*}
  for every $v\in V$ and $\varphi\hspace*{-0.1em}\in\hspace*{-0.1em} C_0^\infty(I)$. Since such a function~is~unique,~${\frac{d\mathbf{u}}{dt}\hspace*{-0.1em}:=\hspace*{-0.1em}\mathbf{w}}$ is  well-defined. By
  \begin{align*}
    \mathbfcal{W}:=W^{1,p,p'}(I,V,V^*)
 :=\big\{\mathbf{u}\in \mathbfcal{X}
    \mid \exists \, \tfrac{d\mathbf{u}}{dt}\in
   L^{p'}(I,V^*)\big\} ,
  \end{align*}
  we denote the \textbf{Bochner--Sobolev space} with respect to the
  evolution triple $(V,H,\textrm{id})$.
\end{definition}

In the context of evolutionary problems, the following
generalized~notions of pseudo-monotonicity and coercivity
(cf.~Definition~\ref{2.1}) are particularly relevant and useful.

\begin{definition}[Bochner pseudo-monotonicity]\label{3.2}
  \hspace*{-0.2em}An
  operator
  ${\mathbfcal{A}\hspace*{-0.15em}:\hspace*{-0.2em}\mathbfcal{X}\hspace*{-0.15em}\cap\hspace*{-0.15em}\mathbfcal{Y}\hspace*{-0.2em}\rightarrow\hspace*{-0.2em}\mathbfcal{X}^*}$ 
  is said to be \textbf{Bochner} \textbf{pseudo-monotone} if for a sequence
  ${(\mathbf{u}_n)_{n\in\mathbb{N}}\hspace*{-0.15em}\subseteq\hspace*{-0.15em}
  \mathbfcal{X}\hspace*{-0.15em}\cap\hspace*{-0.15em}\mathbfcal{Y}}$ from
  \begin{alignat*}{2}
    \mathbf{u}_n\overset{n\rightarrow\infty}{\rightharpoonup}
    &\mathbf{u}\quad &&\text{ in }\mathbfcal{X}.
    \\\mathbf{u}_n\;\;\overset{\ast}{\rightharpoondown}\;\;
    &\mathbf{u}&&\text{ in }\mathbfcal{Y} \quad
    (n\rightarrow\infty),
    \\\mathbf{u}_n(t)\overset{n\rightarrow\infty}{\rightharpoonup}
    &\mathbf{u}(t)\quad &&\text{ in }H\quad\text{for
      a.e. }t\in I,
  \end{alignat*}
  and
  \begin{align*}
    \limsup_{n\rightarrow\infty}{\langle \mathbfcal{A} \mathbf{u}_n,
    \mathbf{u}_n-\mathbf{u}\rangle_{\mathbfcal{X}}}
    \leq 0,
  \end{align*}
  it follows that
  ${\langle
  \mathbfcal{A}\mathbf{u},\mathbf{u}-\mathbf{v}\rangle_{\mathbfcal{X}}\leq
  \liminf_{n\rightarrow\infty}{\langle
    \mathbfcal{A}\mathbf{u}_n,\mathbf{u}_n-\mathbf{v}\rangle_{\mathbfcal{X}}}}$~for~every~${\mathbf{v}\in \mathbfcal{X}}$.
\end{definition}
\begin{definition}[Bochner coercivity]\label{3.20}
  An operator
  ${\mathbfcal{A}: \mathbfcal{X}\cap\mathbfcal{Y}
    \rightarrow\mathbfcal{X}^*}$~is~called:
	\begin{enumerate}
        \item[(i)] \textbf{Bochner coercive with respect to
            $\mathbf{f}\in\mathbfcal{X}^*$ and $u_0\in H$}~if~there is a constant
          $M:=M(\mathbf{f},u_0,\mathbfcal{A})>0$ such that for
          every
          ${\mathbf{u}\in\mathbfcal{X}\cap\mathbfcal{Y}}$~from
		\begin{align*}
		\tfrac{1}{2}\|\mathbf{u}(t)\|_H^2
		+\langle\mathbfcal{A}\mathbf{u}
			-\mathbf{f},\mathbf{u}\chi_{\left[0,t\right]}\rangle_{\mathbfcal{X}}
		\leq \tfrac{1}{2}\|u_0\|_H^2\quad\text{ for a.e. }t\in I,
		\end{align*}
		it follows that $\|\mathbf{u}\|_{\mathbfcal{X}\cap\mathbfcal{Y}}=\|\mathbf{u}\|_{\mathbfcal{X}}+\|\mathbf{u}\|_{\mathbfcal{Y}}\leq M$.
		\item[(ii)] \textbf{Bochner coercive} if it is Bochner coercive with
		respect to $\mathbf{f}$~and~$u_0$, for every
		$\mathbf{f}\in\mathbfcal{X}^*$ and $u_0\in H$.
	\end{enumerate}
\end{definition}

The critical role of the above definitions is that they identify a
vast class of problems for which existence can be established. 
In fact, if  ${\mathbfcal{A}: \mathbfcal{X}\cap\mathbfcal{Y} \to\mathbfcal{X}^*}$ is bounded, Bochner
pseudo-monotone, and Bochner coercive, then the corresponding
evolution problem $\frac{d\bu}{dt}+\mathbfcal{A}\bu=\bff$
is solvable for any initial~datum $u_{0}\in H$. This result was recently obtained in \mbox{\cite[\hspace*{-0.1em}Thm.~\hspace*{-0.1em}4.1]{KR2021}}. 
%

This result is particularly relevant since the difficulty is then
shifted to the verification of the properties of induced operators,
which can be performed time-by-time in the known steady setting. We will not
describe the full result, but we propose a particular, 
simplified setting~enough~to~solve~\eqref{eq:Baldwin--Lomax-alpha}. 

The existence result is  mainly based on the following
proposition giving sufficient conditions which have to be
checked at any fixed~time~slice~${t\in I}$~and which is a particular case  of \cite[Prop.~3.13]{KR2021}.

\begin{proposition}\label{3.1}
  Let $A:V \to V^*$ be an operator. Assume that there exists a number
  $p \in (1,\infty)$ and constants $c_0,c_1>0$ such that\footnote{
    For a pseudo-monotone operator $A:X\to X^*$
  (local) boundedness implies demi-continuity, i.e., $x_n\to x$ in $X$ $(n\to \infty)$ implies $Ax_n\rightharpoonup Ax$ in $X^*$ $(n\to \infty)$,
  hence we do not need here to make any further assumptions of demi-continuity.}:
  \begin{description}
  \item[\textbf{(C.1)}]\hypertarget{C.2} For every $v\in V$ there
    holds 
      \begin{align*}
        \left\|Av\right\|_{V^*}\leq
        c_{0}\|v\|_{V}^{p-1}.
      \end{align*}
    \item[\textbf{(C.2)}]
      \hypertarget{C.3}{}$A:V \rightarrow V^*$ is
      pseudo-monotone.
      \item[\textbf{(C.3)}] \hypertarget{C.4}{} For every $v\in V$
        there holds
      \begin{align*}
      	\langle Av,v\rangle_{V}\ge
      	c_{1}\|v\|_V^p.
      \end{align*}
    \end{description}
    Then, the induced operator
    $\mathbfcal{A}:\mathbfcal{X}\cap \mathbfcal{Y}\rightarrow\mathbfcal{X}^*$, for all $\mathbf{u}\in \mathbfcal{X}\cap\mathbfcal{Y}$ and $ \mathbf{v}\in \mathbfcal{X}$~defined~by 
    \begin{align*}
    	\langle \mathbfcal{A}\mathbf{u},\mathbf{v}\rangle_{\mathbfcal{X}}:=\int_I{\langle A(\mathbf{u}(t)),\mathbf{v}(t)\rangle_{V}\,dt},
    \end{align*}
	is well-defined, bounded, Bochner pseudo-monotone, and Bochner coercive. 
\end{proposition}

On the basis of Proposition~\ref{3.1}, we immediately obtain the following~existence result, which will be used to study the families of
rotational models~just~checking that the conditions
\textup{(\hyperlink{C.1}{C.1})--(\hyperlink{C.3}{C.3})} are satisfied,
after a proper choice of the functional setting.

\begin{theorem}\label{4.2}
  Let 
  $A:V\to V^*$ be an operator satisfying
  \textup{(\hyperlink{C.1}{C.1})--(\hyperlink{C.3}{C.3})}. Then, for arbitrary
  $u_0\in H$ and ${\mathbf{f}\in  L^{p'}(I,V^*)}$, there
  exists a solution $\mathbf{u}\in \mathbfcal{W}$ of the evolution
  equation \begin{align}\label{eq:1}
  	\int_I\Big\langle\frac{d\mathbf{u}}{dt}(t)+A(\mathbf{u}(t)),\mathbf{v}(t)\Big\rangle_{V}&= \int_I{\langle\mathbf{f}(t),\mathbf{v}(t)\rangle_{V}\,dt}\quad\forall\mathbf{v}\in  \mathbfcal{X},\notag\\      
  \mathbf{u}_c(0)&=u_0\quad\text{
  		in }H.\notag
  \end{align}
  Here, the initial condition  has to be understood in the
  sense of the unique continuous representation $\mathbf{u}_c\in C^0(\overline{I},H)$ of $\mathbf{u}\in \mathbfcal{W}$ (cf.~\cite[Prop. 23.23]{zei-IIA}).
\end{theorem}

\section{Weighted spaces}
\label{sec:weighted}
Since \eqref{eq:Baldwin--Lomax-alpha} is a boundary value problem with
the principal part given by a space dependent (and degenerate
at the boundary) operator, a natural functional setting would be that
of weighted Sobolev spaces. Apart from classical Lebesgue and Sobolev
spaces, we will use their weighted counterparts. We follow the notation
from the classical book of Kufner et al.~\cite{Kuf1985}.

A weight  $\varrho$ on $\setR^n$ is a locally integrable function satisfying
almost everywhere $0<\varrho(\bfx)<\infty$. The weighted space
$L^{p}(\Omega,\varrho)$, $1<p<\infty$, is defined as follows
\begin{equation*}
  L^{p}(\Omega,\varrho):=\Big \{\bff:\ \Omega\to\R^{n} \text{ measurable }\fdg
    \int_{\Omega}|\bff(\bfx)|^{p}\, \varrho(\bfx)\,\mathrm{d}\bfx<\infty\Big
    \} .
\end{equation*} 
For $p>1$ we have by using H\"older's inequality that 
\begin{equation*}
  \varrho^{-1/(p-1)}\in L^{1}_{\textup{loc}}(\Omega)\quad\Rightarrow\quad  L^{p}(\Omega,\varrho)\subset
  L^{1}_{\textup{loc}}(\Omega)\subset \mathcal{D}'(\Omega),
\end{equation*}
allowing to work in the standard setting of distributions. It turns
out that $C^{\infty}_{0}(\Omega)$ is dense in $L^{p}(\Omega,\varrho)$
if the weight satisfies $\varrho^{\frac {-1}{p-1}}\in
L^{1}_{\textup{loc}}(\setR^n)$,~see~\cite{Kuf1985}.  In addition,
$L^{p}(\Omega,\varrho)$ is a Banach space when equipped with the norm
\begin{equation*}
  \|\bff\|_{p,\varrho}:=\bigg(\int_{\Omega}|\bff(\bfx)|^{p}\varrho(\bfx)\,\mathrm{d}\bfx\bigg)^{1/p}. 
\end{equation*}
Next, we define weighted Sobolev spaces
\begin{equation*}
  W^{k,p}(\Omega,\varrho):=\left\{\bff:\ \Omega\to\R^{n}\fdg D^{\alpha}\bff\in
    L^{p}(\Omega,\varrho)\text{ for all }\alpha \text{ s.t. } |\alpha|\leq k\right\} ,
\end{equation*} 
equipped with the norm 
\begin{equation*}
  \|\bff\|_{k,p,\varrho}:=\Bigg(\sum_{|\alpha|\leq k}\|D^{\alpha}\bff\|_{p,\varrho}^{p}\Bigg)^{1/p},
\end{equation*}
and, as usual, we define $W^{k,p}_{0}(\Omega,\varrho)$ as follows
\begin{equation*}
  W^{k,p}_{0}(\Omega,\varrho):=\overline{\left\{\bfphi\in
      C^{\infty}_{0}(\Omega)\right\}}^{\|\,.\,\|_{k,p,\varrho}}. 
\end{equation*}

\medskip

In our application the weight $\varrho(\bfx)$ will be a power of the distance
$d(\bfx)\geq0$ of the point $\bx\in \Omega$ from the boundary
$\partial\Omega$. Consequently, we specialize to this setting and
give specific notions regarding these so-called \textit{power-type
  weights}, see Kufner~\cite{Kuf1985}. First, it turns out that
$W^{k,p}(\Omega,d^{\alpha})$ is a separable Banach space provided
$\alpha\in\R$, $k\in\mathbb{N}$ and $1\leq p<\infty$. In this special setting,
since $d(\bx)\geq C_{K}>0$ for each compact $K\subset\subset \Omega$, 
several results are stronger or more precise due to the inclusion
$L^{p}(\Omega,d^{\alpha})\subset L^{p}_{\textup{loc}}(\Omega)$, valid  for all
$\alpha\in \R$.

We recall the following classical result about the distance function
(cf.~\cite{Kuf1985}). 
\begin{lemma}
\label{lem:distance}
Let $\Omega$ be a domain of class $C^{0,1}$, which means that in a
small enough 
neighborhood $\Omega_{P}$, for $P\in \partial\Omega$, 
the boundary $\partial\Omega\cap \Omega_{P}$ can be expressed (after a
rigid rotation)  as $x_{3}=a(x_{1},x_{2})$ for Lipschitz continuous $a$. 
Then, there exist constants
  $0<c_{0},c_{1}\in\R$ such that 
  \begin{equation*}
    c_{0}\,d(\bfx)\leq|a(x')-x_{3}|\leq c_{1}\,d(\bfx)\qquad      \forall\,\bfx=(x',x_{3})\in
    \Omega_{P}. 
  \end{equation*}
\end{lemma}
One of the most relevant properties of the distance function is that
the following embedding holds true
\begin{equation}
\label{eq:embedding-L1}L^{p}(\Omega,d^{\alpha})\subset L^{1}(\Omega)\quad\text{if}\quad
\alpha<p-1.
\end{equation}
It follows directly from H\"older's inequality 
\begin{equation*}
  \int_{\Omega}|f
  |\,\mathrm{d}\bfx=\int_{\Omega}d^{\alpha/p}|f
  |d^{-\alpha/p}d\bfx
  \leq\Big(\int_{\Omega}d^{\alpha}|f|^{p}d\bfx\Big)^{1/p} \Big(\int_{\Omega}d^{-\alpha
    p'/p}d\bfx\Big)^{1/p'},   
\end{equation*}
and using  Lemma~\ref{lem:distance}  the latter integral is finite if and only if 
\begin{equation*}
  \frac{ \alpha\,p'}{p}=\frac{\alpha}{p-1}<1.
\end{equation*}
In the same way we have also that
\begin{equation}
\label{eq:holder}
\forall\,\alpha\in[0,p-1[\qquad  L^{p}(\Omega,d^{\alpha})\subset L^{q}(\Omega)\qquad \forall \,
q\in\big[1,\frac{p}{1+\alpha}\big[.
\end{equation}
As in~\cite[Prop.~9.10]{Kuf1985} it can be shown that:
\begin{lemma}\label{lem:poincare}
  The quantity
  $\Big(\int_{\Omega}d^{\alpha}|\nabla
  \bff|^{p}\,\mathrm{d}\bfx\Big)^{\frac{1}{p}}$ is an equivalent norm
  in $W^{1,p}_{0}(\Omega,d^{\alpha})$, provided that
  $0\leq \alpha<p-1$.
\end{lemma}
In this case
functions from $W^{1,p}_{0}(\Omega,d^{\alpha})$ are zero on
$\partial\Omega$ in the sense that they can be approximated by smooth
functions with compact support.
%
In the sequel we will use certain Hardy--Sobolev
inequalities. Note that inequalities of this kind, when $d$ is replaced
by $|\bx|=d(\bx,0)$ are known as Caffarelli--Kohn--Nirenberg
inequalities~\cite{CKN1984}. 
\begin{lemma}
Let $\Omega \subseteq \setR^n$ be a bounded Lipschitz domain. For $p
\in [1,n)$, $\alpha\neq p-1$ and $q \in [p,\frac {np}{n-p}\big]$ there
exists a constant $c>0$ such that for all $f \in
W^{1,p}_0(\Omega,d^\alpha)$ there holds
\begin{equation}
    \label{eq:dyda}
    \left(\int_{\Omega}d^{\frac qp (n-p+\alpha)-n}|f|^{q}\,d\bfx\right)^{\frac 1q}\leq c\, \left(
      \int_{\Omega}d^{\alpha}|\nabla f|^{p}\,d\bfx\right)^{\frac 1p}. 
  \end{equation}
\end{lemma}
\begin{proof}
  This follows from the definition of the space
  $W^{1,p}_0(\Omega,d^\alpha)$, \cite[Theorem 2.1]{LV16} and the
  classical $(p,\alpha)$ Hardy inequality
\begin{equation}
    \label{eq:Hardy}
    \left(\int_{\Omega}d^{\alpha-p}|f|^{p}\,d\bfx\right)^{\frac 1p}\leq c\, \left(
      \int_{\Omega}d^{\alpha}|\nabla f|^{p}\,d\bfx\right)^{\frac 1p},
  \end{equation}
which is valid for all $p\in(1,\infty)$ and $\alpha\not=p-1$, for
functions in   $W^{1,p}_0(\Omega,d^\alpha)$
  (cf.~\cite{Nec62},
  \cite[Theorem 8.10.14]{kfj}).
\end{proof}
In addition to~\eqref{eq:embedding-L1} and its role in Hardy-type 
inequalities, the critical nature of the power $\alpha=p-1$  also occurs
in the notion of Muckenhoupt weights and their relation with the
maximal function.
\begin{definition}\label{def:muc}
  We say that a weight $\varrho \in L^{1}_{\textup{loc}}(\R^{3})$ belongs to the Muckenhoupt
  class $A_{p}$, for $1<p<\infty$, if there exists $C$ such that
  \begin{equation*}
    \sup_{Q\subset \setR^{n}}\Bigg(
    \fint_{Q}\varrho  (\bfx)\,\mathrm{d}\bfx\Bigg)\Bigg(\fint_{Q}
    \varrho  (\bfx)^{1/(1-p)}\,\mathrm{d}\bfx\Bigg)^{p-1}\leq C,
  \end{equation*}
where $Q$ denotes a cube in $\R^{3}$.
\end{definition}
  The powers of the distance function belong to the class $A_{p}$
according to the following well-known result for general domains (say
it is enough that $\partial \Omega$ is a $n-1$-dimensional closed set,
see \cite{DILTV2019}). Here and in the sequel the boundary will be at
least locally Lipschitz to have the outward unit vector properly
defined.
\begin{lemma}\label{lem: Muckenhoupt}
  The function $\varrho  (\bfx)=\big(d(\bfx)\big)^{\alpha}$ is a Muckenhoupt weight of class
  $A_{p}$ if and only if $-1<\alpha<p-1$.
\end{lemma}
%
%
\subsection{Solenoidal spaces}
A standard approach in fluid mechanics, is to incorporate the
divergence-free constraint directly in the function spaces. These
spaces are built upon completing the space of solenoidal smooth vector
fields with compact support, denoted as
$\bfphi\in C^{\infty}_{0,\sigma}(\Omega)$. For $\alpha\in \R$ define
\begin{equation*}
  \begin{aligned}
   & L^{p}_{\sigma}(\Omega,d^{\alpha}):=\overline{\left\{\bfphi\in
        C^{\infty}_{0,\sigma}(\Omega)\right\}}^{\|\,.\,\|_{p,d^{\alpha}}},
    \\
    &W^{1,p}_{0,\sigma}(\Omega,d^{\alpha}):=\overline{\left\{\bfphi\in
        C^{\infty}_{0,\sigma}(\Omega)\right\}}^{\|\,.\,\|_{1,p,d^{\alpha}}}\,.
  \end{aligned}
\end{equation*}
For $\alpha=0$ they reduce to the classical spaces $L^{p}_{\sigma}(\Omega)$ and
$W^{1,p}_{0,\sigma}(\Omega)$.  Next, we will extensively use the
following extension of classical inequalities linking curl/divergence
and full gradient estimates (cf.~\cite{BB2020}).
\begin{lemma}
  \label{lem:positivity}
  Let $1<p<\infty$ and assume that the weight $\varrho $ belongs to
  the class $A_p$. Then, there exists a constant $C$, depending on the
  domain $\Omega$ and on the weight $\varrho \in A_{p}$, such that
\begin{equation*}
  \|\nabla \bfu\|_{p, \varrho }\leq C ( \|\Div \bfu\|_{p, \varrho }+
  \|\curl \bfu\|_{p, \varrho })\qquad \forall\,\bfu\in W^{1,p}_0(\Omega, \varrho  ).
\end{equation*}
\end{lemma}
In particular, we will use the latter result in the following special form
\begin{corollary}
  For $-1<\alpha<p-1$ there exists a constant $C=C(\Omega,\alpha,p)$ such that
\begin{equation}
  \label{eq:grad-curl-weighted}
  \int_\Omega d^\alpha|\nabla\bfv|^p\,\mathrm{d}\bfx\leq C  \int_\Omega
  d^\alpha|\curl\bfv|^p\,\mathrm{d}\bfx\qquad\forall\, \bfv\in W^{1,p}_{0,\sigma}(\Omega,d^\alpha). 
\end{equation}
\end{corollary}

\section{Application to the rotational turbulence models: the proof of
Theorem~\ref{thm:main-theorem}}\label{sec:sec5}
In this section we verify that the initial boundary value
problem~\eqref{eq:Baldwin--Lomax-alpha}, after a proper selection of
parameters,  and definition of
both the operators and functional spaces, can be put in the
framework of the abstract Theorem~\ref{4.2}. This will be enough to
give a proof of the main result of this paper, that is the existence
of weak solutions in Theorem~\ref{thm:main-theorem}.

In our setting the choice of the natural spaces is determined by the
problem itself which yields, by the a priori estimate obtained by
testing with the velocity $\overline \bfv$, that the integral 
$$\int_{0}^{T}\int_{\Omega}d^{\alpha}|\curl\overline\bfv|^{3}\,d\bx dt$$ is
finite. Hence, for almost all $t\in [0,T]$ the integral
$\int_{\Omega}d^{\alpha}|\curl \overline \bfv|^{3}\,d\bx$ will be finite,
determines the choice for the Banach space $V$.

In order to identify the evolution triple to be used for the proper
formulation, we need to clarify the relationship with the
$L^{2}(\Omega)$ norm. 
We have the following result which immediately derives 
from the basic results on weighted spaces of the previous section.
\begin{lemma}
  \label{lem:Gelfand}
 Let  $\bu\in C^{\infty}_{0,\sigma}(\Omega)$ and 
 $\alpha \in [0,2)$. Then, there exists $C=C(\alpha,\Omega)$ such that
 \begin{equation}
\left(   \int_{\Omega}|\bu|^{2}\,d\bx\right)^{\smash{1/2}} \leq C
\left(\int_{\Omega}d^{\alpha}|\curl\bfu|^{3}\,d\bx\right)^{\smash{1/3}}. 
  \end{equation}
\end{lemma}
\begin{proof}
For $\alpha<2$, combining \eqref{eq:dyda} with
$q=\frac{3p}{3-p+\alpha}$ and
\eqref{eq:grad-curl-weighted},  it follows for~every $p \in (\alpha+1,3)$
\begin{equation*}
  \int_\Omega{|\bfu|^{\frac{3p}{3-p+\alpha}}\,\mathrm{d}\bfx}\leq c\,
  \left( \int_\Omega
    d^\alpha|\nabla\bfu|^p\,\mathrm{d}\bfx\right)^{\frac qp}\leq c
  \left( \int_\Omega 
  d^\alpha|\curl\bfu|^p\,\mathrm{d}\bfx\right)^{\frac qp}, 
\end{equation*}
for all $\bfv\hspace*{-0.05em}\in
\hspace*{-0.05em}C^{\infty}_{0,\sigma}(\Omega) $. Since $2\le
\frac{3p}{3-p+\alpha}$ the assertion follows from H\"older's
inequality as $\Omega$ is bounded.
\end{proof}
Lemma~\ref{lem:Gelfand} shows that one can
work with the following~evolution~triple for all $\alpha\in [0,2)$
\begin{equation*}
(V,H,\textrm{id}):=\big(W^{1,3}_{0,\sigma}(\Omega,d^\alpha),
L^{2}_{\sigma}(\Omega),
\textrm{id}\big). 
\end{equation*}
%
and as functional setting for \eqref{eq:Baldwin--Lomax-alpha}
we use the  following spaces and operators, where  $0\leq \alpha<2$
\begin{equation*}
  \begin{aligned}
    V&:=W^{1,3}_{0,\sigma}(\Omega,d^{\alpha})\qquad
    \|\bfv\|_{V}:=\left(\int_{\Omega}d^{\alpha}|\curl
      \bfv|^{3}\,d\bfx\right)^{1/3}
    \\
    H&:=L^{2}_{\sigma}(\Omega)\qquad
    \|\bfv\|_{H}:=\left(\int_{\Omega}|\bfv|^{2}\,d\bfx\right)^{1/2}
   \\
   \mathbfcal{X}&:=L^3(I,V),\qquad 
   \mathbfcal{Y}:=L^\infty(I,H)
   \\
    \mathbfcal{W}&:=\Big\{\mathbf{u}\in L^{3}(I,V)
    \fdg \exists \, \frac{d\mathbf{u}}{dt}\in
   L^{3/2}(I,V^*)\Big\} ,
\end{aligned}
\end{equation*}
  and define the operator $A:=S+B:V\to V^*$ via 
\begin{equation*}
  \begin{aligned}
    \langle S\bfv, \bfw\rangle_{V}&:=\int_{\Omega} d^\alpha|\curl\bfv|\curl\bfv\cdot \curl\bfw
    \,d\bfx,
    \\
    \langle B \bfv,\bfw\rangle_{V}&:= \int_\Omega
    (\curl \bfv\times\bfv)\cdot\bfw\, d\bfx. 
  \end{aligned}
\end{equation*}
%
The induced operator
$ \mathbfcal{S}:\ \mathbfcal{X}\cap
\mathbfcal{Y}\rightarrow\mathbfcal{X}^*$ inherits the properties of
the operator $S$ (cf.~\cite[Chapter 30]{Zei90B}). Note that $S$ is a
strictly monotone, bounded, coercive, and continuous operator.  These
properties are practically the same known for the $p$-Laplace
operator. In fact, from the definition, one obtains directly the
following two inequalities:
  \begin{equation*}
    \begin{aligned}
      \left\|S\bfv\right\|_{V^*}&\leq
      \|\bfv\|_{V}^{2}\qquad \forall \,\bfv\in V,
      \\
            \langle S\bfv,\bfv\rangle_V&=\|\bfv\|_V^{3}\qquad \forall \,\bfv\in V.
    \end{aligned}
\end{equation*}
The monotonicity of $S $
  derives from the following lemma 
  (cf.~\cite[Lemma~3.3]{BB2020}).
\begin{lemma}
  For smooth enough vector field $\bfomega_{i}$ (it is actually enough that
  $d^\frac{\alpha}{p}\bfomega_i\in L^{p}(\Omega)$, with $1<p<\infty$)
  and for $\alpha\in\R^+$ it holds that
    \begin{equation*}
      \int_\Omega
      (d^{\alpha}|\bfomega_1|^{p-2}\bfomega_1-d^{\alpha}|\bfomega_2|^{p-2}\bfomega_2)
\cdot(\bfomega_1-\bfomega_2)\,\mathrm{d}\bfx\geq0,
\end{equation*}
for any (not necessarily the distance) bounded function such that
$d:\Omega\to\R^+$ for a.e. $\bfx\in\Omega$.
\end{lemma}    
The proof of the above lemma is based on the observation that it can
be proved that 
$d^{\alpha}(|\bfomega_1|^{p-2}\bfomega_1-|\bfomega_2|^{p-2}\bfomega_2)
\cdot(\bfomega_1-\bfomega_2)\geq0$ point-wise. Then weighted
integrability of the functions this used to
prove that the integral is finite.

To treat the operator $B$, and the induced one
$\mathbfcal{B}:\, \mathbfcal{X}\cap
\mathbfcal{Y}\rightarrow\mathbfcal{X}^*$, we need to properly adapt
the estimates on the convective term in weighted spaces and this is
mainly based on the previously Hardy-type
inequalities~\eqref{eq:dyda}.

%
\begin{lemma}[Boundedness of $B$]\label{lem:bd}
  For all $\alpha\in [0, 2)$ the operator ${B:V\to V^{*}}$ is bounded.
  It satisfies
  $\langle B\bfu,\bfv\rangle_{V}\leq c \|\bfu\|_{V}^{2}\|\bfv\|_{V}$ and
  $\langle B\bfu,\bfu\rangle_{V}=0$, for all $\bfu,\bfv\in V$.
\end{lemma}
\begin{proof}
  The proof is based on the estimation of the space integral, by using
  appropriate weighted version of classical Sobolev spaces tools. We
  have in fact, for all smooth functions with compact support the
  following inequality (obtained multiplying and dividing
  a.e. $\x\in\Omega$ by
  the positive function $d^{\alpha/3}$)
  \begin{equation*}
    \begin{aligned}
      &\left|
        \int_{\Omega}( \curl\bfv\times \bfu) \cdot\bfw\,d\x\right|\leq
      \int_{\Omega} d^{-\alpha/6}|\bfu|\, d^{\alpha/3}|\curl
      \bfv|\,d^{-\alpha/6}|\bfw|\,d\x
      \\
      &\quad\leq\left(\int_{\Omega}d^{-\alpha/2}|\bfu|^{3}\,d\x\right)^{1/3}
      \left(\int_{\Omega}d^{\alpha}|\curl
        \bfv|^{3}\,d\x\right)^{1/3}
      \left(\int_{\Omega}d^{-\alpha/2}|\bfw|^{3}\,d\x\right)^{1/3}.
    \end{aligned}
  \end{equation*}
  It remains to show that all $\bfu \in V$ also belong to the weighted
  space $L^3(\Omega, d^{-\frac \alpha 2})$, with a continuous
  embedding. From~\eqref{eq:dyda} it follows for all $p \in [1,3)$,
  $\alpha \in [0,2)$ and
  $\bfu \in W^{1,p}_{0,\sigma}(\Omega,d^{\alpha})$ that
  \begin{equation*}
    \left(\int_{\Omega}d^{- \frac \alpha 2}|\bfu|^{\frac
        {p(6-\alpha)}{2(3-p+\alpha)} }\,d\bfx\right)^{\frac 1q}\leq
    c\, \left( 
      \int_{\Omega}d^{\alpha}|\nabla \bfu|^{p}\,d\bfx\right)^{\frac 1p}, 
  \end{equation*}
  with
  $$
    q:=\frac {p(6-\alpha)}{2(3-p+\alpha)}  <p^{*}.
  $$
  One easily checks that for all $\alpha \in [0,2)$ there exists a $p
  \in (1+\alpha,3)$ such that $3<q  <p^{*}$. Since $\Omega$ is bounded we deduce from this  $V \vnor L^3(\Omega,
  d^{-\frac \alpha 2})$ by using H\"older's inequality.

  Once the
  integral $\int_{\Omega}(\bfu\times \curl\bfv) \cdot\bfw\,d\x$ 
  is
  well-defined for
  $\bfu,\bfv,\bfw\in V$,  
  it
  immediately follows that $\langle B\bfu,\bfu\rangle_{V}=0$ for all
  $\bfu\in V$, since  a.e. in $\Omega$ it holds $(\bfv\times \curl\bfv) \cdot\bfv =0$.
\end{proof}
This is enough for what concerns the growth and coercivity. We
need now to show 
compactness for $B$ in order  to prove pseudo-monotonicity.
\begin{lemma}[Compactness of $B$]
  Let $\alpha \in [0,2)$. Then, the weak convergence
  $\bfu_{n}\rightharpoonup \bfu$ in $V$
  implies 
  (up to a sub-sequence) that
  \begin{equation*}
    B\bfu_{n}\to B\bfu \qquad \text{in }V ^*,
  \end{equation*}
  i.e., the operator $B$ is compact.
\end{lemma}
\begin{proof}
  By the boundedness of the weakly converging sequence 
  $(\bfu_{n})_{n\in \mathbb{N}}\subseteq V$ and by \eqref{eq:holder} we get that
  \begin{equation*}
    \|\bfu_{n}\|_{W^{1,r}(\Omega)}\leq C\qquad \forall \,r
   \in\big[1,\frac{3}{1+\alpha}\big[.
  \end{equation*}
  Hence, by the usual (unweighted) compact Sobolev embedding
  $W^{1,r}(\Omega)\hookrightarrow\hookrightarrow L^{\tilde{r}}(\Omega)
  $, valid for all
  $\tilde{r}<(\frac{3}{1+\alpha})^{*}=\frac{3}{\alpha}$ we get also
  that (up to a sub-sequence)  
  \begin{equation*}
  \bfu_{n}\to \bfu\qquad\text{ a.e. and in }L^{\tilde{r}}(\Omega).
\end{equation*}
By using the definition of
$B$, the properties of the
$\curl$ (with summation over repeated indices), and integration by parts, we have that for all $\bfu,\vv\in V$
\begin{equation*}
  \begin{aligned}
    \langle
    B\bfu,\bfv\rangle&=\int_{\Omega}\epsilon_{jkl}\epsilon_{jlm}u_{k}v_{i}\frac{\partial
    u_{m}}{\partial    x_{l}}
\,d\bfx
=\int_{\Omega}(\delta_{kl}\delta_{im}-\delta_{km}\delta_{il})u_{k}v_{i}\frac{\partial
u_{m}}{\partial
    x_{l}}
\,d\bfx
\\
    &=-\int_{\Omega}u_{k}u_{i}\frac{\partial v_{i}}{\partial
      x_{k}}
    \,d\bfx= -\int_{\Omega}(\bfu\otimes\bfu):\nabla\bfv
    \,d\bfx.
  \end{aligned}
\end{equation*}
Hence, we have
  \begin{equation*}
    \begin{aligned}
      \langle B\bfu_{n},\bfv\rangle-\langle B\bfu,\bfv\rangle
      &=     - \int_{\Omega}(\bfu_{n}\otimes \bfu_{n}):\nabla \bfv-(\bfu\otimes
      \bfu):\nabla \bfv \,d\bfx
      \\
      &=
      -\int_{\Omega}\big((\bfu_{n}-\bfu)\otimes \bfu_{n}\big):\nabla
      \bfv+\big(\bfu\otimes 
      (\bfu_{n}-\bfu)\big): \nabla \bfv\, d\bfx.
    \end{aligned}
  \end{equation*}
By H\"older inequality  we get, as in the proof of Lemma~\ref{lem:bd},
  \begin{equation*}
    \begin{aligned}
      &\left| \int_{\Omega}\big((\bfu_{n}-\bfu)\otimes
        \bfu_{n}\big):\nabla \bfv\, d\x\right 
      | 
      \\
&\leq\left(\int_{\Omega}d^{-\alpha/2}|\bfu_{n}-\bfu|^{3}\, d\x \right)^{1/3}\left(\int_{\Omega}d^{\alpha}|\nabla  
        \bfv|^{3}\, d\x \right)^{1/3}\left(\int_{\Omega}d^{-\alpha/2}|\bfu_{n}|^{3}\, d\x \right)^{1/3}
      \\
      &
      \leq\left(\int_{\Omega}d^{-\alpha/2}|\bfu_{n}-\bfu|^{3}\, d\x \right)^{1/3}\|\bfv\|_{V}\|
      \bfu_{n}\|_{V}.
    \end{aligned}
  \end{equation*}
  We now observe that the last two terms are uniformly bounded, while
  \begin{equation*}
    d^{-\alpha/2}|\bfu_{n}-\bfu|^{3}\to0\qquad\text{a.e. }\x\in \Omega.
  \end{equation*}
  Consequently,  to show that the integral vanishes it is enough to prove that
  for some $q>3$ there holds 
  \begin{equation*}
  \|\bfu_{n}-\bfu\|_{L^{q}(\Omega,d^{-\alpha/2})}\leq C
\end{equation*}
uniformly in $n\in
\mathbb{N}$, which permits to apply the Vitali theorem in the weighted space
$L^{3}(\Omega,d^{-\alpha/2})$. However, this was already obtained in
the proof of Lemma \ref{lem:bd}. The other term in the decomposition of
$\langle B\bfu_{n},\bfv\rangle-\langle B\bfu,\bfv\rangle$ can be
treated in the same way. 
\end{proof}


\begin{proof}[Proof of Theorem~\ref{thm:main-theorem}]
  In the previous lemmas we have proved that $A=S+B$ is continuous and
  pseudo-monotone since it is the sum of 
  a monotone continuous and a compact one.  Collecting the estimates
  we have that in particular that the boundedness and coercivity are
  as follows
  \begin{equation*}
    \begin{aligned}
      \left\|A\bfv\right\|_{V^{*}}&\leq c_{0} \|\bfv\|_{V}^{2},
      \\
      \langle A\bfv,\bfv\rangle_{V}&\ge \|\bfv\|_V^{3},
    \end{aligned}
  \end{equation*}
  since $\langle B\bfv,\bfv\rangle_{V}=0$. Hence all hypotheses from $\mathbf{(C.1)}$ to $\mathbf{(C.3)}$ are
  satisfied.

  This shows that the induced operator $\mathbfcal{A}$ is Bochner pseudo-monotone and
  coercive, hence all the hypotheses of the abstract existence
  Theorem~\ref{4.2} are satisfied. This proves the main result of this
  paper, that is the existence of weak solution in
  Theorem~\ref{thm:main-theorem}.
\end{proof}
\subsection{The case $p>3$}
In this section we show that most of the results of the previous
section can be extended (even with easier proofs) to the system with
the following operator  
\begin{equation*}
    \langle S_{p}\bfv, \bfw\rangle_{V}:=\int_{\Omega} d^\alpha|\curl\bfv|^{p-2}\curl\bfv\cdot \curl\bfw
    \,d\bfx\quad \text{ with } p>3,\  0\leq\alpha<p-1,
  \end{equation*}
while the use of the tools typical of pseudo-monotone
  operators fails for $p<3$.
  We can then prove the following result
  \begin{theorem}\label{thm:main-theorem-p}
  Let $p>3$, $\alpha\in [0,p-1)$, $0<T<\infty$,
  $\overline{ \vv_0 } \in L^2_\sigma(\Omega)$, and
  $\bff\in L^{p'}(0,T;(W^{1,p}_0(\Omega,d^\alpha)^*)$. Then, there
  exists a weak solution to the initial boundary value problem
\begin{equation*}
    \begin{aligned} 
      \partial_{t}\vm+\vortm\times\vm +\curl\big(
      d^\alpha|\vortm|^{p-2}\vortm \big) 
      +\nabla \overline{q}&=\bff&&\quad \mbox{in        $(0,T)\times\Omega$,}
      \\
            \vortm&=\curl\vm&&\quad \mbox{in $(0,T)\times\Omega$,}
      \\
      \Div \vm&=0&&\quad \mbox{in $(0,T)\times\Omega$,}
      \\
      \vm&=\mathbf{0}&&\quad \mbox{on $(0,T)\times\partial\Omega$,}
      \\
      \vm(0)&=\overline{ \vv_{0} }&&\quad \mbox{in $\Omega$,}
    \end{aligned}
\end{equation*}
such that
\begin{equation*}
     \vm\in 
    C([0,T];L^2_\sigma(\Omega))\cap
    L^p(0,T;W^{1,p}_{0,\sigma}(\Omega,d^\alpha))
\end{equation*}
    and for all $t\in[0,T]$
\begin{align*}
    &\frac{1}{2}\|\vm(t)\|^2+\int_0^t\int_\Omega C_\alpha
      d^\alpha(\bx)|{\vortm }    (s,\bx)|^p\,d\bx\,ds
  \\
  &=\frac{1}{2}\|\overline{ \vv_0}
    \|^2+\int_0^t\langle\bff(s),\vm(s)\rangle_{W^{1,p}_0(\Omega,d^\alpha)}ds.
\end{align*}
\end{theorem}

The proof of this result is again just a verification that the
hypotheses of the abstract theorem are satisfied, but we will highlight
the critical role of the parameters. 

First note that for $p>3$ and $0\leq\alpha<p-1$ the inclusion
$W^{1,p}_0(\Omega,d^\alpha)\subset L^{2}(\Omega)$ holds true. 
Directly by H\"older's inequality with $\delta=p/2$, $\delta'=p/(p-2)$, and Hardy
inequality~\eqref{eq:Hardy}  we get 
\begin{equation*}
  \begin{aligned}
    \int_{\Omega}|\bu|^{2}\,d\x&=
    \int_{\Omega}d^{\frac{2}{p}(\alpha-p)}|\bu|^{2}d^{\frac{2}{p}(p-\alpha)}\,d\x
    \\
    &\leq \left( \int_{\Omega}d^{\alpha-p}|\bu|^{p}\,d\x\right)^{\frac{2}{p}}
    \left(
      \int_{\Omega}d^{2\frac{p-\alpha}{p-2}}\,d\x\right)^{\frac{p-2}{p}}
    \\
    &\leq C(p,\alpha) \left(
      \int_{\Omega}d^{\alpha}|\nabla\bu|^{p}\,d\x\right)^{\frac{2}{p}},
  \end{aligned}
\end{equation*}
since $p-\alpha>0$. This shows that one can work with the
evolution triple
\begin{equation*}
(V,H,\textrm{id}):=\big(W^{1,p}_{0,\sigma}(\Omega,d^\alpha),
L^{2}_{\sigma}(\Omega),
\textrm{id}\big)\qquad \text{with}\quad p>3,\ 0\leq\alpha<p-1.
\end{equation*} 
The properties of the operator $S_{p}$ are practically the same as those of the
operator $S$, hence one can directly show
  \begin{equation*}
    \begin{aligned}
      \left\|S_p\bfv\right\|_{V^*}&\leq
      \|\bfv\|_{V}^{p-1}\qquad \forall \,\bfv\in V,
      \\
            \langle S_p\bfv,\bfv\rangle_V&=\|\bfv\|_V^{p}\qquad \forall \,\bfv\in V.
    \end{aligned}
\end{equation*}
On the other hand, the properties of $B$ are to be checked. The
operator is the same as before, but the functional setting is different.

To show the boundedness of $B:\,V\to V^{*}$ for $0\leq\alpha<p-1$, we
proceed as in Lemma~\ref{lem:bd} and we get
  \begin{equation*}
    \begin{aligned}
      &\left|
        \int_{\Omega}( \curl\bfv\times \bfu) \cdot\bfw\,d\x\right|
      \leq
      \int_{\Omega} d^{-\alpha/2p}|\bfu|\, d^{\alpha/p}|\curl
      \bfv|\,d^{-\alpha/2p}|\bfw|\,d\x
      \\
      &\leq\left(\int_{\Omega}d^{-\alpha p'/p}|\bfu|^{2p'}\,d\x\right)^{\frac{1}{2p'}}
      \left(\int_{\Omega}d^{\alpha}|\nabla
        \bfv|^{p}\,d\x\right)^{\frac{1}{p}}
      \left(\int_{\Omega}d^{-\alpha p'/p}|\bfw|^{2p'}\,d\x\right)^{\frac{1}{2p'}}.
    \end{aligned}
  \end{equation*}
  Next, observe that $2p'=2p/(p-1)< p$ is satisfied for
  $p>3$. Consequently, 
  in this case 
  we can directly apply H\"older inequality with exponents $\delta=(p-1)/2$ and
  $\delta'=(p-1)/(p-3)$ to bound the first and third integrals as follows:
\begin{equation*}
  \begin{aligned}
    \int_{\Omega}d^{-\alpha
      /(p-1)}|\bfu|^{2p'}\,d\x=\int_{\Omega}d^{(\alpha-p)
      2/(p-1)}|\bfu|^{2p'}d^{(2p-3\alpha) /(p-1)}\,d\x
    \\
    \leq
\left(    \int_{\Omega}d^{\alpha-p}|\bfu|^{p}\,d\x\right)^{\frac{2}{p-1}}
    \left(\int_{\Omega}d^{(2p-3\alpha) /(p-3)}\,d\x\right)^{\frac{p-3}{p-1}}.
  \end{aligned}
\end{equation*}
The first term from the right-hand side is bounded with
$\left(
  \int_{\Omega}d^{\alpha}|\nabla\bfu|^{p}\,d\x\right)^{\frac{2}{p-1}}$
by using~\eqref{eq:Hardy}, 
while the second is finite if
\begin{equation*}
  \frac{2p-3\alpha}{p-3}>-1\quad\Longleftrightarrow\quad \alpha<p-1.
\end{equation*}
This shows that
\begin{equation*}
\langle  B\bu,\bw\rangle\leq C(\Omega,\alpha,p)\|\bu\|_{V}^{2}\|\bw\|_{V},
\end{equation*}
and the compactness of $B$ follows with the same arguments used in the
previous section (almost everywhere convergence and Vitali
theorem). 
\begin{remark}
  The case $1<p<3$ does not fit with the theory for the reasons we now
  explain. The argument with Hardy inequality as in the previous lemma
  requires $p>3$. If we try to apply the same argument used for $p=3$
  with Hardy--Sobolev inequality~\eqref{eq:dyda}, we can write
  \begin{equation*}
    \begin{aligned}
      &
      \left|
        \int_{\Omega}( \curl\bfv\times \bfu) \cdot\bfw\,d\x\right|
      \\
      &\leq\left(\int_{\Omega}d^{-\alpha p'/p}|\bfu|^{2p'}\,d\x\right)^{\frac{1}{2p'}}
      \left(\int_{\Omega}d^{\alpha}|\curl \bfv|^{p}\,d\x\right)^{\frac{1}{p}}
      \left(\int_{\Omega}d^{-\alpha p'/p}|\bfw|^{2p'}\,d\x\right)^{\frac{1}{2p'}},
    \end{aligned}
  \end{equation*}
  and then estimate the first and third integral with \eqref{eq:dyda}
  for $q=2p'<p^{*}$, which holds for $p>\frac{9}{5}$. Hence, to apply
 \eqref{eq:dyda}  the precise exponent will be
  $q=\frac{p}{p-1}\frac{3p-3-\alpha}{3-p+\alpha}$, 
  and since $q\geq 2p'$ this implies that we have to request for 
  \begin{equation*}
    \alpha\leq\frac{5p-9}{3}.
  \end{equation*}
Since we would like to treat cases with $\alpha$ smaller but
``arbitrarily close'' to $p-1$, the inequality
\begin{equation*}
p-1\leq \frac{5p-9}{3},
\end{equation*}
should be correct. On the other hand the latter can be  satisfied only for $p\geq3$. Since we are out of the range of
permitted $p$ this  shows that the estimate can not be used. Being the
inequalities Hardy--Sobolev inequalities sharp,  this proves that
operator $B$ is not bounded for $\frac{9}{5}<p<3$, when $\alpha$ is close to
$p-1$, hence the basic assumptions to use the pseudo-monotone methods
are not satisfied. The existence of weak solutions, if possible, 
should be obtained with different methods and possibly considering different
weak formulations of the problem. \end{remark}

  
\section*{Acknowledgments}
Luigi C. Berselli was partially supported by a grant of the group
GNAMPA of INdAM and by the University of Pisa within
the grant PRA$\_{}2018\_{}52$ UNIPI: ``\textit{Energy and regularity:
  New techniques for classical PDE problems.}''

\def\cprime{$'$} \def\ocirc#1{\ifmmode\setbox0=\hbox{$#1$}\dimen0=\ht0
  \advance\dimen0 by1pt\rlap{\hbox to\wd0{\hss\raise\dimen0
  \hbox{\hskip.2em$\scriptscriptstyle\circ$}\hss}}#1\else {\accent"17 #1}\fi}
  \def\polhk#1{\setbox0=\hbox{#1}{\ooalign{\hidewidth
  \lower1.5ex\hbox{`}\hidewidth\crcr\unhbox0}}}
\providecommand{\bysame}{\leavevmode\hbox to3em{\hrulefill}\thinspace}
\providecommand{\MR}{\relax\ifhmode\unskip\space\fi MR }
\providecommand{\MRhref}[2]{%
  \href{http://www.ams.org/mathscinet-getitem?mr=#1}{#2}
}
\providecommand{\href}[2]{#2}

\end{document}